\documentclass{amsart}
\usepackage{latexsym,amsmath,amsopn,amssymb,amsthm,amsfonts}

\newtheorem*{theorem*}{Theorem}

\newtheorem{prop}{Claim}
\newtheorem{lemma}{Lemma}
\newtheorem*{lemma*}{Lemma}
\newtheorem*{conj*}{Conjecture}

\theoremstyle{definition}
\newtheorem*{example*}{Example}

\theoremstyle{remark}
\newtheorem{remark}{Remark}

\newcommand{\RR}{\mathbb{R}}

\begin{document}
\title{Riemannian $Spin(7)$ holonomy manifold carries  octonionic-K\"{a}hler structure}
\author{Dmitry V. Egorov}
\address{
Ammosov Northeastern federal university, \newline Kulakovsky str. 48, 677000, Yakutsk, Russia}%
\email{egorov.dima@gmail.com}%
\date{}

\thanks{This work was supported in part by Russian Foundation for Basic
Research (grant 09-01-00598-a) and the Council of the Russian
Federation President Grants (projects NSh-7256.2010.1 and
MK-842.2011.1).}

\sloppy

\begin{abstract}
We prove that  Riemannian $Spin(7)$ holonomy  manifolds  carry
octonionic-K\"{a}hler structure.
\end{abstract}

\maketitle

\section{Introduction}

Let us give a brief definition of the octonionic-K\"{a}hler
structure as we understand it. For more rigorous definition see next
section.

Let $(M,g)$ be a smooth Riemannian manifold. Suppose $V$ is a
$7$-dimensional subbundle  of the vector bundle $\mathrm{End}(TM)$
such that a fiber of $V$ through the point is spanned by almost
complex structures $J_\lambda$ at that point.

We impose two constraints on $V$. First there exists a
non-associative product of  almost complex structures. It
corresponds to the octonionic product. Secondly, the following
formula holds:
$$\nabla_g J_\lambda =
\omega_\lambda^\mu J_\mu,$$  where $\omega\in \mathfrak{g}_2\otimes
\Omega^1(M)$. The $\mathfrak{g}_2$ algebra arises naturally, since
$G_2 = \mathrm{Aut}_\mathbb{R}(\mathbb{O})$.

The defined bundle $V$ over $M$ is called an octonionic-K\"{a}hler
structure on manifold $M$ or we say that $M$ is  an
octonionic-K\"{a}hler manifold. We prove the following theorem.
\begin{theorem*}
Let $M$ be a Riemannian $8$-manifold with holonomy group contained
in  $Spin(7)$; then $M$ is the octonionic-K\"{a}hler manifold.
\end{theorem*}

\begin{remark}
The converse statement is proved in the paper \cite{Irish}. Namely,
it is proved that the holonomy group of the octonionic-K\"{a}hler
manifold  in real dimension $8$ is contained in $Spin(7)$.
\end{remark}
\begin{remark}
Since $Sp(1)\cdot Sp(1) = SO(4)$, this theorem is analogous to the
fact that any oriented Riemannian $4$-manifold  carries the
quaternionic-K\"{a}hler structure. In fact, we can say that parallel
vector cross product of rank $3$ is equivalent to the
$\mathbb{A}$-K\"{a}hler structure, where $\mathbb{A} \cong
\mathbb{H},\mathbb{O}$.
\end{remark}

\section{The octonionic-K\"{a}hler structure}
We adopt the definition of the quaternionic-K\"{a}hler manifold
given in \cite[Proposition 14.36]{Besse}.

Suppose $(M,g)$ is a smooth oriented Riemannian manifold. $M$ is
called an {\it octonionic-K\"{a}hler manifold} if there exists an
open cover $U_i$ of $M$ and almost complex structures $J_\lambda$,
$\lambda = 1,\ldots,7$ on each $U_i$ such that
\begin{enumerate}

  \item there exists a non-associative $\times$-product of almost
  complex structures:
  $J_\lambda~\times~J_\mu~=~J_{\lambda\times\mu}$, where
  $\lambda\times\mu$ corresponds to the  product of imaginary unit octonions enumerated by
natural numbers from $1$ to $7$;

  \item $\nabla_g J_\lambda = \omega_\lambda^\mu J_\mu$,
  where $\omega\in \mathfrak{g}_2\otimes \Omega^1(M)$;

  \item for any point  $p\in U_i\cap U_j$
  almost complex structures
  $J_\lambda$, $\lambda = 1,\ldots,7$ span the same vector subspace of the
  $\mathrm{End}(T_pM)$.

  \item metric  $g$ is Hermitian with respect to each $J_\lambda$;

\end{enumerate}

\section{The $Spin(7)$-structure}\label{section_cross}
Let us briefly recall  definition of the $Spin(7)$-structure, for
details see \cite{Joyce}. Define a $3$-form  $\varphi_0$ on $\RR^7$
by
\begin{equation}\label{g2_form}
\begin{array}{c}
\varphi_0 = e^{0145}  + e^{0167}  + e^{2345} + e^{2367} + e^{0246} -
e^{0257} - e^{1346}+ e^{1357}\\- e^{0347}- e^{0356}- e^{1247}-
e^{1256}+ e^{0123}+ e^{4567}.
\end{array}
\end{equation}
By $e^{ijk}$ denote  $e^i\wedge e^j\wedge e^k$, where $e^i$ is the
unit orthogonal coframe. The subgroup of $GL(8,\RR)$ preserving
$\varphi_0$ and orientation is called a $Spin(7)$ group.

Let $M$ be an oriented closed $8$-manifold. Suppose there exists a
global $4$-form $\varphi$ such that pointwise it coincides with
$\varphi_0$; then $M$ is called a $Spin(7)$-manifold or we say that
$M$ carries the $Spin(7)$-structure. Since $Spin(7)\subset SO(8)$,
orientation and the Riemannian metric are uniquely determined by the
$Spin(7)$-structure. If the $Spin(7)$-structure is parallel with
respect to  metric connection of the induced metric, then
$\mathrm{Hol}(M)\subseteq Spin(7)$.

Recall that the first complete $Spin(7)$ holonomy Riemannian metrics
were constructed  in \cite{Bryant}, the first compact in
\cite{Joyce}.

\section{Cross products}
In this section we follow \cite{Gray}. Let $(M,g)$ be a Riemannian
manifold with $Spin(7)$-structure $\varphi$. Suppose a multilinear
alternating smooth map $P:TM\times TM\times TM\to TM$ such that it
is  compatible with metric $g$:
\begin{equation}\label{cross_def1}
g(P(e_1, e_2,e_3),e_i) = 0,\quad  i = 1,2,3;
\end{equation}
\begin{equation}\label{cross_def2}
|P(e_1, e_2,e_3)|^2 = \det g(e_i,e_j),\quad |e|^2 = g(e,e).
\end{equation}
Then $P$ is called a {\it vector cross product}.

The cross product is uniquely determined by the $Spin(7)$-structure
$\varphi$:
\begin{equation}\label{form2prod} \varphi(e_1,e_2,e_3,e_4) = g(
P(e_1, e_2,e_3),e_4),
\end{equation}
where $g$ is induced by $\varphi$. The converse is also true. The
cross product induces the Riemannian metric and the
$Spin(7)$-structure.

The following properties of cross product will be useful for us.

{\bf I.}\quad If cross product is parallel with respect to the
metric connection, then the holonomy group of $M$ is a subgroup of
$Spin(7)$.

{\bf II.}\quad Suppose $u$ and $v$ are fixed vectors; then cross
product determines almost complex structure on the orthogonal
complement to $u$ and $v$ by the following formula:
$$
P(u,v,w) = Jw.
$$

{\bf III.}\quad Composition rule of cross products is described by
the following lemma  \cite[Lemma $4.4.3$]{Karigiannis}.
\begin{lemma*}
Let $(M,g)$ be a Riemannian manifold with  $Spin(7)$-structure
$\varphi$ and cross product $P$; then
\begin{equation}\label{prod_iter}
\begin{array}{l}
P(a,b,P(u,v,w)) = \\
-\: g(a\wedge b, u\wedge v)w - \varphi(a,b,u,v)w + g(b,w)P(a,u,v) -
g(a,w)P(b,u,v)\\
-\: g(a\wedge b, v\wedge w)u - \varphi(a,b,v,w)u + g(b,u)P(a,v,w) -
g(a,u)P(b,v,w)\\
-\: g(a\wedge b, w\wedge u)v - \varphi(a,b,w,u)v + g(b,v)P(a,w,u) -
g(a,v)P(b,w,u).
\end{array}
\end{equation}
Here $ g(a\wedge b, w\wedge u) = g(a,w)g(b,u) - g(a,u)g(b,w)$.

\end{lemma*}

\section{Proof of the main theorem}

Let $e_i$ be a local unit orthogonal frame such that locally
$Spin(7)$-structure $\varphi$ is  of the form \eqref{g2_form}.

Suppose local almost complex structures $J_\lambda$ are determined
by the following formulae:
\begin{equation}\label{almost}
\begin{array}{c}
P(e_0,e_\lambda, v) = J_\lambda v;\\
J_\lambda e_0 = e_\lambda.
\end{array}
\end{equation}
Hereafter Greek indices range over natural numbers from $1$ to $7$.
We also assume that $\lambda$, $\mu$ and $\nu$ are pairwise distinct
if the contrary is not stated.

\begin{prop}There exists a product of local almost complex
structures denoted by $\times$ such that
\begin{equation}\label{times_product} J_\lambda\times J_\mu =
J_{\lambda\times\mu},\end{equation} where $\lambda\times\mu$
corresponds to the product of imaginary unit octonions enumerated by
natural numbers from $1$ to $7$.
\end{prop}
We split the proof into the following lemmata.

\begin{lemma}For any $\lambda$ and $\mu$ there exists $\nu$:
$$ P(e_0, e_\lambda
,e_\mu ) = J_\lambda e_\mu = e_{\nu}.
$$
\end{lemma}
\begin{proof}The proof follows from \eqref{g2_form} and \eqref{form2prod}.
\end{proof}

\begin{lemma}
By definition, put
\begin{equation}\label{dot_product}
e_{\lambda\times\mu} = P(e_0, e_\lambda ,e_\mu);\quad
e_{\lambda\times\lambda} = -e_0; \quad e_{-\lambda} = -e_\lambda.
\end{equation}
Then
\begin{equation}\label{octonion_product}
e_{\lambda\times\mu}  = -\delta_{\lambda\mu}e_0 +
\gamma_{\lambda\mu}^\nu e_\nu,
\end{equation}
where $\delta$ is the Kronecker delta and $\gamma$ are structure
constants of $\mathfrak{g}_2$.
\end{lemma}
\begin{proof}
Cross product $P$ of rank $3$ induces cross product
$P(e_0,\cdot,\cdot)$ of rank $2$ on $e^\perp_0$. Induced product
determines the $G_2$-structure and corresponding structure
constants.
\end{proof}

It is well-known that \eqref{octonion_product} is the product rule
of imaginary unit octonions. Using
\eqref{almost},\eqref{times_product},\eqref{dot_product} and
\eqref{octonion_product}, we define a $\times$-product such that it
determines the  octonionic product of local almost complex
structures.

The Claim $1$ is proved.

\begin{remark}
Generally, the $\times$-product of almost complex structures does
not coincide with their composition, because the latter is
associative.

We conjecture that $\times$-product coincides with the product
defined in \cite{OctonionMatrix}, where the matrix representation of
octonions was constructed.
\end{remark}

\begin{prop}
$\nabla_g J_\nu = \omega_\nu^\mu J_\mu$,   where $\omega\in
\mathfrak{g}_2\otimes \Omega^1(M)$.
\end{prop}
\begin{proof}
Differentiating $P(e_\lambda,e_\mu,v) = J_{\lambda\times\mu}v$, we
have:
\begin{equation}\label{prod_nabla}
P(\nabla_g e_\lambda,e_\mu,v) + P(e_\lambda,\nabla_g e_\mu, v) =
(\nabla_g J_{\lambda\times\mu})v;
\end{equation}
Here we use that $\nabla_g P=0$. There exists a matrix $\rho \in
\mathfrak{so}(8)\otimes \Omega^1(M)$ such that:
\begin{equation}\label{nabla_e}
\nabla_g e_i = \rho^j_i e_j,\quad i = 0,\ldots,7.
\end{equation}
Substituting \eqref{nabla_e} in \eqref{prod_nabla}, we get:
$$
P(\rho^k_\lambda e_k,e_\mu,v) + P(e_\lambda,\rho^l_\mu e_l, v) =
(\nabla_g J_{\lambda\times\mu})v
$$
or $$\rho^k_\lambda J_{k\times\mu} + \rho^l_\mu J_{\lambda\times l}
= \nabla_g J_{\lambda\times\mu}.
$$
Last identity implies that there exists $\omega$ such that $\nabla_g
J_{\lambda\times\mu} = \omega_{\lambda\times\mu}^\mu J_\mu$.
Moreover $\omega$ solves the system of equations:
\begin{equation}\label{slau}
\omega^{\lambda\times\mu}_\lambda + \omega^{\mu\times \lambda}_\mu =
\omega^\mu_{\lambda\times\mu};
\end{equation}
\begin{equation}\label{slau2}
-\omega^\mu_\nu = \omega^\nu_\mu.
\end{equation}
Equations \eqref{slau},\eqref{slau2} form the  $7$-dimensional
representation of $\mathfrak{g}_2$.
\end{proof}

\begin{prop}
For any point  $p\in U_i\cap U_j$ almost complex structures
$J_\lambda$, $\lambda = 1,\ldots,7$ span the same vector subspace of
the   $\mathrm{End}(T_pM)$.\end{prop}
\begin{proof}
The proof follows form the fact that the definition of almost
complex structures \eqref{almost} is linear with respect to the
local frame.
\end{proof}

\begin{prop}
The Riemannian metric $g$ is Hermitian with respect to any of almost
complex structures $J_\lambda$ determined by  \eqref{almost}.
\end{prop}
\begin{proof}
The proof follows from the identity
$$g(e_{\lambda\times\mu},e_{\lambda\times\nu}) =
g(e_\mu,e_\nu).$$ In turn, this identity  follows from
\eqref{prod_iter} and \eqref{dot_product}. This corresponds to the
existence of division in $\mathbb{O}$.
\end{proof}

The proof of the theorem follows from Claims $1$--$4$.

\section{The twistor theory}
The vector bundle $V$ defined in the introduction is associated to
the $6$-sphere bundle $S$, since $G_2/SU(3) = S^6$. The fiber of $S$
through a point of $M$ parametrizes almost complex structures at
that point.

The main theorem of the twistor theory states that total space of
the sphere bundle has integrable complex structure iff base manifold
has restricted curvature \cite{Atiyah,Salamon,Verbit}. If a
$6$-sphere has integrable complex structure, then one can apply the
twistor theory to $Spin(7)$-manifolds.

\end{document}